\documentclass[12pt]{amsart}
\usepackage[dvips]{graphicx}
\usepackage{epsfig}
\usepackage{amssymb}
\usepackage[usenames, dvipsnames]{color}
\usepackage{hyperref}
\usepackage{verbatim}
\usepackage[normalem]{ulem}

\theoremstyle{plain}
\newtheorem{theorem}{Theorem}[section]

\theoremstyle{definition}

\theoremstyle{remark}

\newcommand{\bR}{{\mathbb R}}

\def\e{\epsilon}

\def\p{\partial}

\def\na{\nabla}
\def\al{\alpha}

\def\O{\Omega}

\def\be{\begin{equation}}
\def\ee{\end{equation}}
\def\bes{\begin{equation*}}
\def\ees{\end{equation*}}
\def\bali{\begin{aligned}}
\def\eali{\end{aligned}}
\def\al{\begin{aligned}}
\def\eal{\end{aligned}}

\def\lab{\label}

\def\2O{\underline{\O}}

\numberwithin{equation}{section}

\makeatletter
\def\dashint{\operatorname%
{\,\,\text{\bf--}\kern-.98em\DOTSI\intop\ilimits@\!\!}}
\makeatother

\begin{document}

%\linenumbers

\title[blow up solution]{A blow up solution of the Navier-Stokes equations with  a
critical force}

\thanks{}

\author[Q. S. Zhang]{Qi S. Zhang}

\address[]{Department of mathematics, University of California, Riverside, CA 92521,
USA}

\email{qizhang@math.ucr.edu}



\begin{abstract}
 A forced solution $v$ of the  Navier-Stokes equation in any open domain with no slip
 boundary condition is constructed. The scaling factor of the forcing term is the
 critical order $-2$. The velocity, which is smooth until its final blow up moment,
 is in the energy space through out. Since most physical forces from a point source
 in nature are regarded as order $-2$, such as Coulomb force, Yukawa force, this
 result indicates possible singularity formation under these kind of forces. The result
  even holds for some log subcritical forces or some forces in the standard critical space $L^\infty_t
 L^{3/2}_x$, including the explicit force:  $F=- \delta \frac{e^{-|x|^2}}{(|x|^2 + T-t) \,[1+
| \ln (|x|^2 + T-t)|]} (1, 0, 0) $ for any small $\delta>0$. The result can also be considered as a step in Scheffer's plan proposed in \cite{Sc} p49.

\end{abstract}
%\today
\maketitle

%\tableofcontents

\section{Statement of result}

The  Navier-Stokes equations (NS) modeling the motion of viscous incompressible
fluids in a domain $D \subset \bR^3$ are
\be
\lab{nse}
\mu \Delta v -  v \nabla v - \nabla P -\partial_t v =F, \quad div \, v=0, \quad
\text{in} \quad D \times (0, \infty),  \quad v(\cdot, 0)=v_0(\cdot)
\ee The unknowns are the velocity $v$ and the pressure $P$.
 $\mu>0$ is the viscosity constant, taken as $1$ here; $F$ is a given forcing term.
 Suitable boundary conditions are imposed if $\partial D \neq \emptyset$.

  Due to Leray \cite{Le2}, if $D=\bR^3$, $v_0 \in L^2(\bR^3)$, $F=0$,  the Cauchy
  problem has a solution in the energy space $(v, \nabla v) \in   (L^\infty_t L^2_x,
  L^2_{x t})$. However, in general we do not know if such solutions stay bounded or
  regular for all $t>0$ given regular initial values. The uniqueness of such
  solutions is also unknown so far.

 Solutions beyond the energy space or super critical forcing terms have been
 investigated before. Recall that for the forcing term, the space $L^\infty_t
 L^{3/2}_x$ is critical. Let us recall that a forcing term in local $L^q_tL^p_x$
 space is super critical, critical and sub critical if $\frac{3}{p}+ \frac{2}{q} >,
 =, < 2$ respectively. In general, under the standard scaling, a forcing term  is
 super critical, critical and sub critical if its scaling order is less than, equal
 to or greater to $-2$ respectively. A simpler analogy is the forced heat equation $\Delta v - \partial_t v =F$.

 For a forcing term in the supercritical space $L^1_t L^2_x$, in the paper
 \cite{ABC22},   non-uniqueness of \eqref{nse} in the energy space is established.
 In the homogeneous case $F=0$, non-uniqueness was proven for some solutions in the
 space $L^\infty_t L^2_x$ in \cite{BV}. Recently, an instantaneous  finite time blow
 up example was constructed in \cite{LPYZZZ} for a special cusp domain. The solution
 is in the energy space and the forcing term (Ampere force) is actually subcritical.
 These results indicate that if pertinent function spaces or the domain are
 sufficiently singular, then something singular can happen to the solutions. In a
 related development, forced blow up solutions for some hypo-dissipative Navier
 Stokes equations was studied in \cite{CMZ}.

In a previous note \cite{Zq23}, a forced solution $v$ of the axially symmetric
Navier-Stokes equation in a finite cylinder $D$ with suitable boundary condition is
constructed. The divergence free, rotational forcing term has an order of scaling
which is a little worse than the critical order $-2$. See also \cite{BY} for an extension. One can argue that this is not
too surprising  since for any smooth, divergence free vector field $v$ that blows up
in a given final time, the following trivial equation holds before blow up:
\[
\Delta v -\partial_t v -v \nabla v = F  \equiv \Delta v -\partial_t v -v \nabla v
\] Then this vector field $v$ is a solution to a forced Navier-Stokes equation.  In
this way one can produce blow up solutions with slightly super critical forcing term
in the localized space $L^\infty_t L^{3^-/2}_x$. Here $3^-$ is any positive number
strictly less than $3$.
An example from \cite{Zq23} is a follows. Let $\phi=\phi(y)$ be a smooth, divergence
free vector field, which is compactly supported in $D$. For example $\phi =
\text{curl} \, X$ where $X$ is any smooth vector field compactly supported in $D$.
Define
$
v=v(x, t) = \ln (1/(1-t)) \phi (x/\sqrt{1-t}).
$ Then $v$ is a blow up solution of the above Navier-Stokes equation with $0$
pressure and $0$ boundary value in $D \times [0, 1)$.
 It is easy to check that $v$ is in the energy space and $F \in L^\infty_t
 L^{3^-/2}_x([0, 1] \times D)$. The problem with this almost effortless process is
 that the forcing term $F$ is super critical or nonphysical, as pointed out long time ago on p50 of \cite{Sc}. In that paper, V. Scheffer constructed a blow up solution of the forced Navier Stokes equations in the energy space with a very singular force $F$.
 %The main novelty is that $F \cdot v \ge 0$ in our sign convention so the force actually slows down the growth of the velocity. His result is a step in his plan in constructing singular solutions of the Navier Stokes equations in the energy space, but the force is considered unnatural and bizarre c.f. line 5, p50 \cite{Sc}.
 The main novelty is that $F \cdot v \ge 0$ in our sign convention so the force might seem to slow down the growth of the velocity. Scheffer's result is a step in his plan to construct singular solutions of the Navier-Stokes equation in the energy space, but he himself characterizes the force as "a bizarre function that would not come up in real life situations", c.f. line 5, p. 50 in \cite{Sc}.

 In this paper, we improve the results in \cite{Zq23} and, to some degree, \cite{Sc} by showing some critical forces
 or even log sub-critical forces can also generate finite time singularity, although we  should add that the sign of $F \cdot v$ is not controlled. Since
 most physical forces from point sources in nature are regarded as order $-2$,  this
 result indicates possible singularity formation under these natural forces.
The first main result has two parts. The force in part (a) has a principal
part modeled on physical forces which have a scaling factor of $-2$. The force in
part (b) lives in the mathematically more convenient critical space $L^\infty_t
L^{3/2}_x$. We mention that the result does not and can not directly address the
regularity problem which requires the forcing term to be subcritical or even smooth.
Although there may be some relation between the two due to the fact that those forces
which cause the blow up are so close to being genuinely sub-critical such as
$L^\infty_t L^{3^+/2}_x$ functions. The result can also be considered as a step in Scheffer's plan.

\begin{theorem}
\lab{prblows}
Let $D$ be any open domain in $\mathbb{R}^3$ containing the origin $0$, and $Q$ be
the space time domain $D \times [0, T)$, for any given $T>0$.

(a). There exists a forcing term $F$ whose principal term locally equals
$
 -\frac{e^{-|x|^2}}{|x|^2 + T-t} (1, 0, 0)
$ such that the Navier Stokes equations \eqref{nse} have a smooth, compactly
supported solution $v$ in the energy space in $Q$, which blows up at time $T$.

(b). There exists a forcing term $F \in L^\infty_t L^{3/2}_x(\mathbb{R}^3 \times [0,
T))$  whose principal term locally equals $-\frac{e^{-|x|^2}}{(|x|^2 + T-t) \,[1+
| \ln (|x|^2 + T-t)|]} (1, 0, 0)$ such that the Navier Stokes equations \eqref{nse} have a smooth solution $v$ in
the energy space in $\mathbb{R}^3 \times [0, T)$, which blows up at time $T$.
\end{theorem}

Let us outline the idea of the proof. We first construct a blow up solution $h$ of a vector valued heat equation with a critical forcing term, with a log blow up rate. Then we will transform $h$ to a divergence free vector field $v$ by the  formula $curl curl (\Delta^{-1}) h$. The error term turns out to be a gradient field which absorbs the most unnatural term into the pressure.
We will show that the nonlinear term $v \nabla v$ scales log worse than $-1$ and turn it into lower order part of the force.  Note that we will use the pressure in a substantial way but we have not used the nonlinearity fully. So future improvement may be possible.

Even though the principal part of the force in the theorem is explicit, one can still argue that the total force is somewhat artificial since the lower order terms are not explicit.  However, with a little further effort, aided by a previous result \cite{Zq04}, we can prove the existence of  blow up solutions with explicit forces and zero initial value as follows.

\begin{theorem}
\lab{prblows2}
 There exists a sufficiently small number $\delta_0>0$ such that for all $\delta \in (0, \delta_0]$ the following statement holds.  The Navier Stokes equations \eqref{nse} with the $L^\infty_t L^{3/2}_x$ force $F=- \delta \frac{e^{-|x|^2}}{(|x|^2 + T-t) \,[1+
| \ln (|x|^2 + T-t)|]} (1, 0, 0) $ and  $0$ initial value have a smooth solution $v$ in
the energy space in $\mathbb{R}^3 \times [0, T)$, which blows up at time $T$.
\end{theorem}

Note that the magnitude of this force $F$ is already less singular than natural forces from point sources.
This theorem is  built on the previous theorem by a perturbation argument using the Stokes kernel. In the proof of the  previous theorem, we actually showed that the linear Stokes equation with the given force in $L^\infty_t L^{3/2}_x$  has a solution $Z$ that blows up in finite time. Now we look for a blow up solution of the Navier Stokes equations in a neighborhood of $Z$ under a suitable norm. It turns out that we just need to solve a perturbed Navier Stokes equations whose extra terms are subcritical. So a small data bounded classical solution exists. Combining with $Z$, the claimed blow up solution is constructed.

\section{ Proof}

\proof (of Theorem \ref{prblows}.)

 Without loss of generality we take $T=1$. We will also use $c,  C$ to denote positive generic constants which may change value.

(a). By scaling we can assume the ball $B(0, 2) \subset D$.

{\it Step 1.}
We start by solving the following vector valued heat equation in the whole space.
\be
\begin{cases}
\lab{eqh100}
\Delta h(x, t) - \partial_t h(x, t) = f \equiv -\frac{e^{-|x|^2}}{|x|^2 + 1-t} (1, 0,
0), \quad t \in [0, 1)\\
h(\cdot, 0)=0.
\end{cases}
\ee By the standard heat kernel formula in $\mathbb{R}^3$, we know that $h$ is given
by
\[
h=(h_1, 0, 0)
\]where $h_1$ is given by
\be
\lab{hfml}
\al
h_1 (x, t) &=\frac{1}{(4 \pi)^{3/2}} \int^t_0 \int \frac{1}{(t-s)^{3/2}}
e^{-\frac{|x-y|^2}{4(t-s)}} \frac{e^{-|y|^2}}{|y|^2 + 1-s} dyds\\
&=\frac{1}{(4 \pi)^{3/2}} \int^t_0 \int \frac{1}{s^{3/2}} e^{-\frac{|x-y|^2}{4 s}}
\frac{e^{-|y|^2}}{|y|^2 + 1-t+s} dyds,\\
&=\frac{1}{(4 \pi)^{3/2}} \int^t_0 \int  e^{-\frac{|(x/\sqrt{s})-y|^2}{4}} \frac{e^{-
s|y|^2}}{ s |y|^2 + 1-t+s} dyds
\eal
\ee where  changes of variable $s$ to $t-s$ first and then $y$ to $y \sqrt{s}$  are
made.
This formula implies, since $ t \in [0, 1)$ that
\be
\lab{h0log}
\al
h_1(0, t) & \ge \frac{1}{(4 \pi)^{3/2}}  \int  e^{-1.25 |y|^2} \int^t_0\frac{1}{ s
(|y|^2+1) + 1-t} dsdy\\
&=\frac{1}{(4 \pi)^{3/2}}  \int  \frac{e^{-1.25 |y|^2}}{1+|y|^2}  \ln
\frac{t(1+|y|^2)+1-t}{1-t} dy \\
& \ge \frac{1}{(4 \pi)^{3/2}}  \int  \frac{e^{-1.25 |y|^2}}{1+|y|^2} dy \,   \ln
\frac{t}{1-t}.
\eal
\ee Therefore $h_1(0, t)$ blows up like $\ln[1/(1-t)]$ as $t \to 1^-$, i.e. for some
$c>0$,
\be
\lab{h0ln1-t}
h_1(0, t) \ge c |\ln(1-t)|.
\ee  On the other hand, notice that for $t \in [0, 1]$, we have
\[
\al
&\int^t_0  \frac{1}{(t-s)^{3/2}} e^{-\frac{|x-y|^2}{4(t-s)}}ds = \int^t_0
\frac{1}{s^{3/2}} e^{-\frac{|x-y|^2}{4s}}ds\\
&= \frac{1}{|x-y|} \int^{t/|x-y|^2}_0  \frac{1}{l^{3/2}} e^{-\frac{1}{4l}}dl
\le \frac{1}{|x-y|} \int^{1/|x-y|^2}_0  \frac{1}{l^{3/2}} e^{-\frac{1}{4l}}dl \le
\frac{C e^{-|x-y|^2/5 }}{|x-y|}
\eal
\]
It is then quickly seen from \eqref{hfml} that
\be
\lab{h1sj}
\al
h_1(x, t) &\le C \int \frac{e^{-|x-y|^2/5}}{|x-y|} \frac{e^{-|y|^2}}{|y|^2+ 1-t }
dy\\
&\le C \int_{\{|x-y| \le |x|/2\} } ... dy + C \int_{\{|x|/2 \le |x-y| \le 3 |x|\}}
... dy +
C \int_{\{|x-y| \ge 3|x|\}} ... dy\\
&\le C e^{-c |x|^2} [1+|\ln (|x|^2+1-t)|].
\eal
\ee

{\it Step 2.}

Let $\phi=\phi(x, t)$ be a smooth cut-off function such that $\phi=1$ in $B(0, 1)
\times [1/2, 1]$ and $\phi=0$ outside of $B(0, 2) \times [0, 1]$.
We further look for  solutions of \eqref{nse} in the form of
\be
\lab{solnsf}
\al
v&=v(x, t) \equiv \frac{1}{4 \pi} curl
\left( \phi(x, t)\,  curl \int \frac{1}{|x-y|} h(y, t) dy \right)\\
&\equiv curl
\left( \phi(x, t)\,  curl ( (-\Delta)^{-1} h)  \right).
\eal
\ee Due to the decay of $h_1$ near infinity, given in \eqref{h1sj}, we know
$(-\Delta)^{-1} h$ is well defined and so is $v$. By vector calculus
\be
\lab{v=cc}
\al
v&= \phi \, curl \, curl ((-\Delta)^{-1} h) + \nabla \phi \times curl \,(
(-\Delta)^{-1} h)\\
&=\phi h + \phi \nabla ( div \, (-\Delta)^{-1} h) + \nabla \phi \times curl \,(
(-\Delta)^{-1} h).
\eal
\ee The main point is that here the second term on the right hand side is a gradient
vector field in $B(0, 1) \times [1/2, 1)$ and the third term is bounded and smooth on
the closed time interval $[0, 1]$. On the semi-open time interval $[0, 1)$, we have
\[
\al
\Delta v -\p_t v&= (\Delta h -\p_t h) \phi + 2 \nabla \phi \nabla h + h (\Delta \phi
-\p_t \phi)\\
&\qquad +\nabla \left[ (\Delta -\p_t) ( div \, (-\Delta)^{-1} h)\right] +
(\Delta -\p_t)\left[(\phi-1) \nabla ( div \, (-\Delta)^{-1} h) \right]\\
&\qquad +(\Delta -\p_t)\left[\nabla \phi \times curl \,( (-\Delta)^{-1} h)\right]
\eal
\]Writing
\be
\lab{fbded}
\al
f_b&=2 \nabla \phi \nabla h + h (\Delta \phi -\p_t \phi)+(\Delta -\p_t)\left[(\phi-1)
\nabla ( div \, (-\Delta)^{-1} h) \right]\\
&\qquad
 +(\Delta -\p_t)\left[\nabla \phi \times curl \,( (-\Delta)^{-1} h)\right],
\eal
\ee
\be
\lab{pres}
P=(\Delta -\p_t) ( div \, (-\Delta)^{-1} h),
\ee and using \eqref{eqh100}, we find that
\[
\Delta v -\p_t v=\nabla P + f \phi +f_b.
\]Note that $f_b$ is smooth and bounded through out since its support is outside the
unit parabolic cube $Q_1 =B(0, 1) \times [1/2, 1]$ where all functions involved are
smooth and bounded. Hence $v$ is a solution of the following forced Navier-Stokes
equations with no-slip boundary condition:
\be
\lab{nseF}
 \Delta v -  v \nabla v - \nabla P -\partial_t v =F, \quad div \, v=0, \quad
 \text{in} \quad D \times [0, 1),  \quad v(\cdot, 0)=0,
\ee where
\be
\lab{forceF}
F=f \phi -v \nabla v + f_b.
\ee

To finish the proof we need to check the following assertions.

1. $v$ becomes unbounded as $t \to 1^-$ but smooth for $t<1$.

2. $v$ is in the energy space in the closed time interval $[0, 1]$.

3.
 $v \nabla v$ scales logarithmically worse than $-1$. So $f$, whose scaling order is
 $-2$, is the dominant force since, as mentioned, $f_b$ is bounded and smooth.

 Let us prove assertion 1 first.
 It it suffices to show that $v(0, t)$ becomes unbounded as $t \to 1^-$. So we can
 just concentrate on the domain where $\phi=1$. According to \eqref{v=cc}, in this
 domain we have
 \[
 \al
 v&=  h +  \nabla ( div \, (-\Delta)^{-1} h)\\
 &=(h_1, 0, 0) + \nabla ( div \, ((-\Delta)^{-1} h_1, 0, 0))\\
 &=(h_1, 0, 0) + \nabla ( \p_1(-\Delta)^{-1} h_1)\\
 &=((-\Delta) (-\Delta)^{-1} h_1, 0, 0) +  ( \p^2_1(-\Delta)^{-1} h_1, \p_2
 \p_1(-\Delta)^{-1} h_1, \p_3 \p_1(-\Delta)^{-1} h_1).
 \eal
 \]Therefore
 \be
 \lab{v=cc2}
 v=( -(\p^2_2+\p^2_3) (-\Delta)^{-1} h_1, \, \p_2 \p_1(-\Delta)^{-1} h_1, \, \p_3
 \p_1(-\Delta)^{-1} h_1).
 \ee We claim that the function $(\p^2_2+\p^2_3) (-\Delta)^{-1} h_1$ is unbounded as
 $t \to 1^-$. If it were bounded, then $(\p^2_1+\p^2_2) (-\Delta)^{-1} h_1$ and
 $(\p^2_1+\p^2_3) (-\Delta)^{-1} h_1$ would be bounded too by radial symmetry of
 $h_1$ and change of variables $(x_1, x_2, x_3) \to (x_3, x_1, x_2)$ and $(x_1, x_2,
 x_3) \to (x_2, x_1, x_3)$. For example, making the change of variables $z_1=x_2$,
 $z_2=x_1$ and $z_3=x_3$, we have
 \[
 \al
 &(\p^2_2+\p^2_3) (-\Delta)^{-1} h_1(x, t)=\frac{1}{4\pi}(\p^2_{x_2}+\p^2_{x_3}) \int
 \frac{ h_1(y, t)}{\sqrt{|x_1-y_1|^2+|x_2-y_2|^2+|x_3-y_3|^2}} dy\\
 &=\frac{1}{4\pi}(\p^2_{z_1}+\p^2_{z_3}) \int \frac{ h_1(y_1, y_2, y_3,
 t)}{\sqrt{|z_2-y_1|^2+|z_1-y_2|^2+|z_3-y_3|^2}} dy\\
 &=\frac{1}{4\pi}(\p^2_{z_1}+\p^2_{z_3}) \int \frac{ h_1(y_2, y_1, y_3,
 t)}{\sqrt{|z_2-y_2|^2+|z_1-y_1|^2+|z_3-y_3|^2}} dy, \quad (y_2 \to y_1, y_1 \to
 y_2),\\
 &=\frac{1}{4\pi}(\p^2_{z_1}+\p^2_{z_3}) \int \frac{ h_1(y_1, y_2, y_3,
 t)}{\sqrt{|z_1-y_1|^2+|z_2-y_2|^2+|z_3-y_3|^2}} dy, \quad (h_1 \, \text{ is
 radial}),\\
 &=\frac{1}{4\pi}(\p^2_{x_1}+\p^2_{x_3}) \int \frac{ h_1(y,
 t)}{\sqrt{|x_1-y_1|^2+|x_2-y_2|^2+|x_3-y_3|^2}} dy.  \quad (z_i \to x_i, i=1, 2,
 3.)
 \eal
 \]Notice that due to the change of variables, the above argument does not infer that the three functions are equal pointwise. What being shown is that
 \[
 \Vert (\p^2_2+\p^2_3) (-\Delta)^{-1} h_1 \Vert_\infty = \Vert (\p^2_1+\p^2_3) (-\Delta)^{-1} h_1 \Vert_\infty= \Vert (\p^2_1+\p^2_2) (-\Delta)^{-1} h_1 \Vert_\infty.
 \]
  Adding these three functions together would yield
 \[
 2(\p^2_1+ \p^2_2+\p^2_3) (-\Delta)^{-1} h_1=- 2 h_1
 \]is bounded. This is false since $h_1$ is unbounded by \eqref{h0ln1-t}. Therefore
 the claim is true, implying that $v$ is unbounded as $t \to 1^-$. The smoothness of
 $v$ when $t<1$ follows from \eqref{solnsf} and the smoothness and exponential decay
 of $h_1$ when $t<1$, which follows from \eqref{hfml} via routine calculations. So
 assertion 1 is proven.

Next we prove assertion 2. i.e. for all $t \in [0, 1)$, the following holds.
\be
\lab{enerv}
\int_{D} |v|^2(x, t) dx +
\int^t_0 \int_{D} |\na v|^2(x, s) dx ds \le C.
\ee Here $C$ is a constant independent of $t$. According to \eqref{v=cc2}, at each
time level $t$,  $v$ is a vector field whose components are the Riesz transforms of
$h_1$, the solution of \eqref{eqh100}. By the standard property of the Riesz
transforms, it suffices to prove the energy bound for $h_1$, namely
\be
\lab{enerh}
\int h_1^2(x, t) dx +
\int^t_0 \int |\na h_1|^2(x, s) dx ds \le C.
\ee
From \eqref{h1sj}, it clear that $\int_{D} h_1^2(x, t) dx \le C.$
Taking the gradient on the first equality of
\eqref{hfml}, we deduce
\be
\al
|\nabla h_1 (x, t)| & \le C \int^t_0 \int \frac{1}{(t-s)^{2}}
e^{-\frac{|x-y|^2}{4.5(t-s)}} \frac{e^{-|y|^2}}{|y|^2 + 1-s} dyds\\
&=C \int^t_0 \int \frac{1}{s^2} e^{-\frac{|x-y|^2}{4.5 s}} \frac{e^{-|y|^2}}{|y|^2 +
1-t+s} dyds,\\
&\le C \int \int^t_0 \frac{1}{s^2} e^{-\frac{|x-y|^2}{4.5 s}} ds
\frac{e^{-|y|^2}}{|y|^2 + 1-t} dy\\
&\le C \int  \frac{e^{-\frac{|x-y|^2}{5}}}{|x-y|^2}  \frac{e^{-|y|^2}}{|y|^2 + 1-t}
dy\\
&\le C \int_{\{|x-y| \le |x|/2\} } ... dy + C \int_{\{|x|/2 \le |x-y| \le 3 |x|\}}
... dy +
C \int_{\{|x-y| \ge 3|x|\}} ... dy\\
&\le C e^{-c |x|^2} [1+\frac{1}{|x|+\sqrt{1-t}}].
\eal
\ee This implies that
\[
\int^1_0 \int |\na h_1|^2(x, t) dx dt \le C,
\]proving \eqref{enerh} and hence \eqref{enerv}. Thus assertion 2 is proven.

Now we prove assertion 3.

By \eqref{v=cc2}, we see that
\be
\lab{dv1=}
\nabla v_1 = -(R_2 R_2 + R_3 R_3) \nabla h_1, \quad \nabla v_2 = R_2 R_1  \nabla
h_1,
\quad \nabla v_3 = R_3 R_1  \nabla h_1
\ee Here $R_i$ are the Riesz integrals and $R_iR_j$ are Riesz transforms. Using \eqref{hfml}, we find
\[
\p^2_{x_i x_j} h_1 (x, t) =\frac{1}{(4 \pi)^{3/2}} \int^t_0 \int
\frac{1}{(t-s)^{3/2}} \p_{x_i} e^{-\frac{|x-y|^2}{4(t-s)}} \p_{y_j}
\frac{e^{-|y|^2}}{|y|^2 + 1-s} dyds,
\]which infers
\be
\lab{hfmd2}
\al
&\left| \p^2_{x_i x_j} h_1 (x, t) \right|
\\
&\le C \int^t_0 \int \frac{1}{s^2} e^{-\frac{|x-y|^2}{4.5 s}-|y|^2}
\left[\frac{1}{(|y|^2 + 1-t+s)^{3/2}} +1 \right] dyds\\
& \le C \int  \frac{e^{-\frac{|x-y|^2}{5}-|y|^2}}{|x-y|^2} \left[\frac{1}{(|y|^2 +
1-t)^{3/2}} +1 \right] dy\\
&\le \frac{C e^{-c |x|^2}}{ (|x|^2+1-t)^{3/2}} \int_{\{|x-y| \le |x|/2\} }
\frac{1}{|x-y|^2} dy + \frac{C e^{-c |x|^2}}{|x|^2} \int_{ |x|/2 < |x-y| \le 3 |x|\}}
\frac{1}{(|y|^2+1-t)^{3/2}} dy \\
&\qquad +
C e^{-c |x|^2} \int_{\{|x-y| \ge 3|x|\}} \frac{e^{-c |y|^2}}{|y|^2
(|y|^2+1-t)^{3/2}}\,  dy.
\eal
\ee Hence
\[
\left| \p^2_{x_i x_j} h_1 (x, t) \right|
\le C \frac{e^{-c |x|^2}}{ |x|^2+1-t} + C \frac{e^{-c |x|^2}}{ |x|^2} \ln\left(1+
\frac{|x|^3}{(1-t)^{3/2}} \right).
\]From here, by considering the separate cases when $|x|^2 \le 1-t$ and $|x|^2 >1-t$,
we infer
\[
\left| \p^2_{x_i x_j} h_1 (x, t) \right|
\le C \frac{e^{-c |x|^2}}{ |x|^2+1-t}\left[1+\ln\left(1+ \frac{|x|^3}{(1-t)^{3/2}}
\right) \right].
\]

 Using this and integration by parts, we obtain, after splitting the domains of
 integration as before, that
\[
\al
&|R_i R_j \nabla h_1(x, t) | \le C \int \frac{1}{|x-y|^2} |\nabla^2 h_1(y, t)| dy + C |\nabla h_1(x, t)| \\
&\le C \int \frac{1}{|x-y|^2} \frac{e^{-c |y|^2}}{ |y|^2+1-t} \left[1+\ln\left(1+
\frac{|y|^3}{(1-t)^{3/2}} \right) \right]dy +  C |\nabla h_1(x, t)| \\
&\le \frac{C}{1+ |x|^2} [1+\frac{1}{|x|+\sqrt{1-t}}] \left[1+\ln\left(1+
\frac{|x|^3}{(1-t)^{3/2}} \right) \right].
\eal
\]This can be made rigorous by doing the usual cut-off for the singular kernel $\partial_{x_i}\partial_{x_j} \frac{1}{|x-y|}$ in the small balls $B(x, \e)$ and let $\e$ go to $0$ after integration by parts. By \eqref{dv1=}, this infers
\be
\lab{dvjie}
|\nabla v(x, t)| \le \frac{C}{1+ |x|^2} [1+\frac{1}{|x|+\sqrt{1-t}}]
\left[1+\ln\left(1+ \frac{|x|^3}{(1-t)^{3/2}} \right) \right].
\ee Using \eqref{v=cc2} and integration by parts again, we have, similar to
\eqref{h1sj}, that
\[
\al
|R_i R_j  h_1(x, t) | &\le C \int \frac{1}{|x-y|^2} |\nabla h_1(y, t)| dy +C |h_1(x, t)|\\
&\le C \int \frac{1}{|x-y|^2} \frac{e^{-c |y|^2}}{ |y|+\sqrt{1-t}} dy +C |h_1(x, t)|\\
&\le \frac{C}{1+ |x|^2} [1+|\ln(|x|^2+1-t)|],
\eal
\]which infers, via \eqref{v=cc2} once more,
\be
\lab{vjie}
|v(x, t)| \le \frac{C}{1+ |x|^2} [1+|\ln(|x|^2+1-t)|].
\ee

Therefore $v \nabla v$ scales logarithmic worse than $-1$, proving assertion 3.
This completes the proof of part (a).
\medskip

(b). The proof is similar to (a). The difference is that we choose a forcing term
that is logarithmic weaker than that in part (a) to make it living in $L^\infty_t
L^{3/2}_x$.

This time we solve the following vector valued heat equation in the whole space.
\be
\begin{cases}
\lab{eqh101}
\Delta h(x, t) - \partial_t h(x, t) = f \equiv -\frac{e^{-|x|^2}}{(|x|^2 + 1-t) \,[1+
| \ln (|x|^2 + 1-t)|]} (1, 0, 0), \quad t \in [0, 1)\\
h(\cdot, 0)=0.
\end{cases}
\ee
Then $h$ is given by
\[
h=(h_1, 0, 0),
\]
\be
\lab{hfml2}
\al
h_1 (x, t) &=\frac{1}{(4 \pi)^{3/2}} \int^t_0 \int \frac{1}{(t-s)^{3/2}}
e^{-\frac{|x-y|^2}{4(t-s)}} \frac{e^{-|y|^2}}{(|y|^2 + 1-s) \,[1+ | \ln (|y|^2 +
1-s)|]} dyds\\
&=\frac{1}{(4 \pi)^{3/2}} \int^t_0 \int  e^{-\frac{|(x/\sqrt{s})-y|^2}{4}} \frac{e^{-
s|y|^2}}{ [s (|y|^2+1) + 1-t]\,[1+ | \ln (s(|y|^2+1) + 1-t)|]} dyds
\eal
\ee where similar changes of variable as in part (a)  are made.
This formula implies, for $ t \in [0.8, 1)$ that
\be
\lab{h0log}
\al
h_1(0, t) & \ge \frac{1}{(4 \pi)^{3/2}}  \int  e^{-1.25 |y|^2}
\int^t_0\frac{1}{[s (|y|^2+1) + 1-t]\,[1+ | \ln (s(|y|^2+1) + 1-t)|]} dsdy\\
& \ge C \int_{0.5 \le |y| \le 1}  \frac{e^{-1.25 |y|^2}}{1+|y|^2}  |\ln|
\left[\frac{1+\ln (t(1+|y|^2)+1-t)}{1+|\ln(1-t)|} \right] dy \\
& \ge C    \ln \ln \frac{0.5}{1-t},
\eal
\ee i.e. $h_1(0, t)$ blows up like $\ln \ln[1/(1-t)]$ as $t \to 1^-$,
Now we just define again
\be
\lab{solnsf2}
\al
&v=v(x, t) \equiv curl
\left(  curl ( (-\Delta)^{-1} h)  \right),\\
&P=(\Delta -\p_t) ( div \, (-\Delta)^{-1} h),\\
&F=f -v \nabla v.
\eal
\ee Then $v$ is a solution of \eqref{nse} with $D=\mathbb{R}^3$ and $0$ initial
value.
The rest of the proof follows that of part (a) step by step with only small changes.
First $v$ still blows up at time $t=1$ since the proof only needs the properties that
$h$ blows up at $t=1$ and that $h$ is radially symmetric.  The energy bound for $v$
follows since $h$ and $\nabla h$ have faster decay than that in part (a).  Similarly
$v \nabla v$ decays faster than that in part (a). It is now easy to see the forcing
term $F=f - v \nabla v $ is in the space $L^\infty_t L^{3/2}_x$.  This completes the
proof of part (b) of Theorem \ref{prblows}.
\qed

\proof (of Theorem \ref{prblows2}.)

Again we take $T=1$.
We modify \eqref{eqh101} in the proof of the previous theorem as follows.  Solve the heat equation:
\be
\begin{cases}
\lab{eqh102}
\Delta h(x, t) - \partial_t h(x, t) = F \equiv - \delta \frac{e^{-|x|^2}}{(|x|^2 + 1-t) \,[1+
| \ln (|x|^2 + 1-t)|]} (1, 0, 0), \quad t \in [0, 1),\\
h(\cdot, 0)=0.
\end{cases}
\ee  where $\delta>0$ will be chosen later. Now we just define
\be
\lab{solnsf3}
\al
&Z=Z(x, t) \equiv curl
\left(  curl ( (-\Delta)^{-1} h)  \right),\\
&P=(\Delta -\p_t) ( div \, (-\Delta)^{-1} h).
\eal
\ee Then $Z$ is a solution of the (linear) Stokes problem  in $\mathbb{R}^3 \times [0, 1)$, with $0$ initial
value:
\be
\lab{stkeq}
\Delta Z- \nabla P -\partial_t Z =F, \quad div \, Z=0, \quad
\text{in} \quad \mathbb{R}^3 \times [0, 1),  \quad Z(\cdot, 0)=0
\ee Just like Theorem \ref{prblows2}, $Z$ is in the energy space and it blows up at $t=1$. Moreover, the following estimates hold.
\be
\lab{zjie}
|Z(x, t)| \le  \delta \frac{C}{1+ |x|^2} [1+|\ln(|x|^2+1-t)|].
\ee
\be
\lab{dzjie}
|\nabla Z(x, t)| \le  \delta \frac{C}{1+ |x|^2} \left[1+\frac{1}{|x|+\sqrt{1-t}}\right]
\left[1+\ln\left(1+ \frac{|x|^3}{(1-t)^{3/2}} \right) \right].
\ee  The proofs are just like those of \eqref{vjie} and \eqref{dvjie}. Actually the  $\ln$  terms can be replaced by $\ln \ln$ terms on the right hand side, which will not be written out for simplicity.

Now we look for a solution
\be
\lab{vzu}
v=Z +u
\ee of the Navier Stokes equations
\be
\lab{nser3}
\Delta v  -v \nabla v- \nabla P -\partial_t v =F, \quad div \, v=0, \quad
\text{in} \quad \mathbb{R}^3 \times [0, 1),  \quad v(\cdot, 0)=0.
\ee Direct calculation tells us that $u$ has to satisfy the equations
\be
\lab{ueq}
\begin{cases}
\Delta u  - u \nabla u -Z \nabla u - u \nabla Z- \nabla P -\partial_t u =Z \nabla Z, \quad div \,
 u=0, \quad
\text{in} \quad \mathbb{R}^3 \times [0, 1), \\
 u(\cdot, 0)=0.
 \end{cases}
\ee Here the pressure $P$ has been renamed.

Observe that \eqref{ueq} is a perturbed Navier Stokes equations with a forcing term. The key point is that the perturbation and the force are all subcritical by \eqref{zjie} and \eqref{dzjie} . Indeed, $\nabla Z$ as potential term scales log worse than $-1$ while $-2$ is the critical scale.  Next, $Z$ as the drift, has a log scale, while $-1$ is the critical scale for drift terms.  Finally $Z \nabla Z$ scales log worse than $-1$, again below the critical scale $-2$. This implies that \eqref{ueq} has a bounded regular solution if $Z \nabla Z$ is small in suitable sense.
We start to give a proof of this statement now.

Consider the metric space $X$ for weakly divergence free vector fields $u=u(x, t)$ in $\mathbb{R}^3 \times [0, 1]$ with the weighted norm
\be
\lab{normx}
\Vert u \Vert_X = \sup \{ (1+|x|)^{2} |u(x, t)| \, | \, (x, t) \in  \mathbb{R}^3 \times [0, 1]\}.
\ee  Let $K=K(x, t; y, s)$ be the Stokes kernel. Write $\partial_i K = \partial_{y_i} K(x, t; y, s)$.
Consider the mapping on $X$:
\be
\lab{defmap}
\al
(M u)(x,t)&=  \underbrace{\int^t_0 \int \partial_i K \,  u_i u(y, s) dyds}_{T_1} +
\underbrace{\int^t_0 \int \partial_i K \, Z_i u(y, s) dyds}_{T_2}\\
 &\qquad+ \underbrace{\int^t_0 \int \partial_i K \, u_i Z(y, s) dyds}_{T_3} + \underbrace{\int^t_0 \int \partial_i K \, Z_i Z(y, s) dyds}_{T_4}.
\eal
\ee In the above $u_i$, $Z_i$ are the components of $u$ and $Z$ and summations are taken for repeated indices. From the standard decay properties of $K$ c.f. \cite{So}:
\be
\lab{kdk<}
\al
&|K(x, t; y, s)| \le \frac{C}{(|x-y| + \sqrt{t-s})^3}, \\
&|\nabla K(x, t; y, s)| \le \frac{C}{(|x-y| + \sqrt{t-s})^4 },
\eal
\ee the mapping $M$ is well defined for $u \in X$ and the function $Z$. This has been done for broader function classes in \cite{Zq04}. In fact one only needs $u$ to decay faster than $C/(1+|x|^{1^+})$ where $1^+$ is any number greater than $1$.

We will use the contraction mapping principle to prove that $M$ has a fixed point in the ball
$B(0, \eta) \subset X$ when $\eta$ is a sufficiently small positive number. To this end, we  first need to prove that $M$ maps $B(0, \eta)$ to $B(0, \eta)$ when $\eta$ and $\delta$ are sufficiently small, i.e.
\be
\lab{mbb}
M(B(0, \eta)) \subset B(0, \eta).
\ee We will bound $T_i$, $i=1, ..., 4$ respectively with $T_4$ being the worst. From \eqref{kdk<} and \eqref{zjie}, we see that
\be
\lab{t41}
|T_4| \le \int^t_0 \int |\nabla K | \, |Z|^2 dyds \le
C \delta^2 \int^t_0 \int \frac{(1+|y|^2)^{-2} [|\ln (|y|+\sqrt{1-s})|+1]^2}{(|x-y| + \sqrt{t-s})^4 } dyds
\ee For $0<s< t \in (0, 1]$, we have
\[
|\ln (|y|+\sqrt{1-s})|+1 \le C \left( |\ln |y||+1 \right),
\]and
\be
\lab{tintk}
\al
\int^t_0 &\frac{[|\ln (|y|+\sqrt{1-s})|+1]^2}{(|x-y| + \sqrt{t-s})^4 } ds \le
\int^t_0 \frac{C}{(|x-y|^2 + s)^2 } ds \left( |\ln |y||+1 \right)^2\\
&=C \left[\frac{1}{|x-y|^2} - \frac{1}{|x-y|^2+t} \right] \left( |\ln |y||+1 \right)^2\\
& \le C \left[\frac{1}{|x-y|^2} - \frac{1}{|x-y|^2+1} \right] \left( |\ln |y||+1 \right)^2\\
& = \frac{C}{|x-y|^2 (|x-y|^2+1)} \left( |\ln |y||+1 \right)^2.
\eal
\ee Substituting this to \eqref{t41}, we deduce
\be
\lab{t42}
\al
|T_4| &\le C \delta^2 \int \frac{ \left( |\ln |y||+1 \right)^2 (1+|y|^2)^{-2}  }{|x-y|^2 (|x-y|^2+1)}  dy \\
&\equiv C \delta^2 \underbrace{\int_{D_1}  ... dy}_{I_1} +  C \delta^2 \underbrace{\int_{D_2}  ... dy}_{I_2} +  C \delta^2 \underbrace{\int_{D_3}  ... dy}_{I_3},
\eal
\ee where, as before $D_1 =\{ y \, | \, |x-y| \le |x|/2 \}$, $D_2 =\{ y \, | \, |x|/2 <|x-y| \le 3 |x| \}$ and $D_3 =\{ y \, | \, |x-y| > 3 |x| \}$. We will bound $I_1$, $I_2$ and $I_3$ respectively.  First
\be
\lab{i1<}
\al
I_1 &= \int_{D_1} \frac{ \left( |\ln |y||+1 \right)^2 (1+|y|^2)^{-2}  }{|x-y|^2 (|x-y|^2+1)}  dy \le  \frac{C}{(1+|x|^2)^2} \int_{D_1} \frac{ \left( |\ln |x-y||+1 \right)^2  }{|x-y|^2 (|x-y|^2+1)}  dy \\
&=\frac{C}{(1+|x|^2)^2} \int^{|x|/2}_0 \frac{(|\ln \rho|  +1)^2 \rho^2}{\rho^2 (\rho^2+1)} d\rho, \qquad \rho=|x-y|\\
&\le \frac{C}{(1+|x|^2)^2}.
\eal
\ee Next
\be
\lab{i2<}
\al
I_2 &= \int_{D_2} \frac{ \left( |\ln |y||+1 \right)^2 (1+|y|^2)^{-2}  }{|x-y|^2 (|x-y|^2+1)}  dy \le \frac{C}{1+|x|^2} \int_{|y| \le 4 |x|} \frac{ \left( |\ln |y||+1 \right)^2 (1+|y|^2)^{-2} }{|y|^2}  dy\\
&\le \frac{C}{1+|x|^2}.
\eal
\ee Similarly
\be
\lab{i3<}
\al
I_3 &= \int_{D_3} \frac{ \left( |\ln |y||+1 \right)^2 (1+|y|^2)^{-2} }{|x-y|^2 (|x-y|^2+1)}  dy \le \frac{C}{1+|x|^2} \int_{|y| \ge  |x|/2} \frac{ \left( |\ln |y||+1 \right)^2 (1+|y|^2)^{-2} }{|y|^2}  dy\\
&\le \frac{C}{1+|x|^2}.
\eal
\ee Substituting \eqref{i1<}, \eqref{i2<} and \eqref{i3<} into \eqref{t42}, we deduce
\be
\lab{t43}
|T_4| \le \frac{C \delta^2}{1+|x|^2}.
\ee

Next, for $u \in B(0, \eta)$, we deduce, in similar manner, that
\be
\lab{t23}
\al
|T_2| + |T_3| &\le 2 \int^t_0 \int |\nabla K | \, |u| \, |Z| dyds \le
C \delta \eta \int^t_0 \int \frac{(1+|y|^2)^{-1} [|\ln (|y|+\sqrt{1-s})|+1]}{(|x-y| + \sqrt{t-s})^4 (1+|y|^{2})} dyds \\
 &\le \frac{C \delta \eta}{1+|x|^{2}}.
\eal
\ee
For $u \in B(0, \eta)$,  we have, likewise,
\be
\lab{t1}
\al
|T_1|  &\le \int^t_0 \int |\nabla K | \, |u|^2 dyds \le
C  \eta^2  \int^t_0 \int \frac{ 1}{(|x-y| + \sqrt{t-s})^4 (1+|y|^{2})^2} dyds \\
&\le C  \eta^2   \int \frac{C}{|x-y|^2 (|x-y|^2+1)} \frac{ 1}{(1+|y|^{2})^2} dy\\
 &= C  \eta^2   \int_{|x-y|<|x|/2} ...dy + C  \eta^2   \int_{|x-y| \ge |x|/2} ...dy \le \frac{C \eta^2}{1+|x|^{2}}.
\eal
\ee Plugging \eqref{t43}, \eqref{t23} and \eqref{t1} into \eqref{defmap}, we deduce
\[
|M \, u |(x, t) \le \frac{C (\delta+\eta)^2}{1+|x|^{2}}.
\] This proves the containment \eqref{mbb} when we choose $C  (\delta+\eta)^2<\eta$.

Our next task is to prove that the mapping $M$ is a contraction in $B(0, \eta)$. To this end, pick $u^{(1)}, u^{(2)}$ in $B(0, \eta)$. Then, by \eqref{defmap},
\be
\lab{u1-u2}
\al
&(u^{(1)}-u^{(2)})(x, t)\\
&=
\underbrace{\int^t_0 \int \partial_i K \,  (u^{(1)}_i- u^{(2)}_i) u^{(1)}(y, s) dyds}_{T_1} +
\underbrace{\int^t_0 \int \partial_i K \, u^{(2)}_i (u^{(1)}-u^{(2)})(y, s) dyds}_{T_2}\\
 &\qquad + \underbrace{\int^t_0 \int \partial_i K \, Z_i (u^{(1)}-u^{(2)})(y, s) dyds}_{T_3} + \underbrace{\int^t_0 \int \partial_i K \, (u^{(1)}_i-u^{(2)}_i) Z(y, s) dyds}_{T_4}
\eal
\ee
Observe that
\[
\al
|T_1|+|T_2| &\le
2 \Vert u^{(1)}-u^{(2)} \Vert_X \int^t_0 \int |\nabla K| \frac{\eta}{(1+|y|^{2})^2} dyds \\
&\le C \eta \Vert u^{(1)}-u^{(2)} \Vert_X \int^t_0 \int \frac{1}{(|x-y| + \sqrt{t-s})^4 } \frac{1}{(1+|y|^{2})^2} dyds\\
&\le C \eta \Vert u^{(1)}-u^{(2)} \Vert_X  \int \frac{1 }{|x-y|^2 (|x-y|^2+1)} \frac{1}{(1+|y|^{2})^2} dy.
\eal
\]Here we just used the calculation in \eqref{tintk}. By standard estimate of breaking the space into the domains $\{ y \, | \, |x-y| \ge |x|/2 \}$ and its complement, we find that
\be
\lab{2t12}
|T_1|+|T_2| \le C \eta \Vert u^{(1)}-u^{(2)} \Vert_X \frac{1}{1+|x|^{2}}.
\ee Similarly, by \eqref{zjie},
\[
\al
|T_3|+|T_4| &\le 2 \Vert u^{(1)}-u^{(2)} \Vert_X \int^t_0 \int |\nabla K| \, |Z| \frac{1}{(1+|y|^{2})} dyds\\
&\le
 C\Vert u^{(1)}-u^{(2)} \Vert_X \int^t_0 \int |\nabla K| \frac{\delta (1+|y|^2)^{-1} [1+|\ln(|y|^2+1-s)|]}{1+|y|^2} dyds \\
&\le C\Vert u^{(1)}-u^{(2)} \Vert_X \int^t_0 \int |\nabla K| \frac{\delta (1+|y|^2)^{-1}  [1+|\ln|y||]}{1+|y|^{2}} dyds\\
&\le C \delta \Vert u^{(1)}-u^{(2)} \Vert_X \int^t_0 \int \frac{(1+|y|^2)^{-1}  [1+|\ln|y||]}{(|x-y| + \sqrt{t-s})^4 } \frac{1}{1+|y|^{2}} dyds\\
&\le C \delta \Vert u^{(1)}-u^{(2)} \Vert_X  \int \frac{ [1+|\ln|y||] }{|x-y|^2 (|x-y|^2+1)} \frac{1}{(1+|y|^2)^2} dy.
\eal
\]Hence
\be
\lab{2t34}
|T_3|+|T_4| \le C \delta \Vert u^{(1)}-u^{(2)} \Vert_X \frac{1}{|1+|x|^{2}}.
\ee Bring \eqref{2t12} and \eqref{2t34} to \eqref{u1-u2}, we infer
\[
|(u^{1}-u^{2})|(x, t) \le C (\delta+\eta) \Vert u^{(1)}-u^{(2)} \Vert_X \frac{1}{1+|x|^{2}}.
\]Therefore the mapping $M$ is a contraction on $B(0, \eta)$  when we choose $C (\delta+\eta)<1$. Now we can conclude that the mapping has a fixed point $u \in B(0, \eta)$ which satisfies
\be
\lab{intequ}
\al
u(x,t)&=  \int^t_0 \int \partial_i K \,  u_i u(y, s) dyds +
\int^t_0 \int \partial_i K \, Z_i u(y, s) dyds\\
 &\qquad+ \int^t_0 \int \partial_i K \, u_i Z(y, s) dyds + \int^t_0 \int \partial_i K \, Z_i Z(y, s) dyds.
\eal
\ee Since $|u(x, t)|<\eta/(1+|x|^{2})$, we know that $u \in L^\infty_tL^2_x(\mathbb{R}^3 \times [0, 1])$. Due to the boundedness of $u$ and ln type singularity of $Z$, the standard parabolic singular integral theory tells us $|\nabla u| \in L^2_t L^2_x(\mathbb{R}^3 \times [0, 1])$. Indeed the second derivatives of the Stokes kernel is a combination of the second derivatives of the heat kernel and its Riesz transforms, which, regarded as integral kernels, are bounded from space-time $L^2$ to itself. Using this and \eqref{dzjie}, we can carry out integration by parts in the convolutions in \eqref{intequ} to obtain
 \be
\lab{intequ2}
\al
u(x,t)&=  -\int^t_0 \int K \,  u_i  \partial_i u(y, s) dyds -
\int^t_0 \int  K \, Z_i  \partial_i u(y, s) dyds\\
 &\qquad - \int^t_0 \int  K \, u_i  \partial_i Z(y, s) dyds - \int^t_0 \int  K \, Z_i \partial_i Z(y, s) dyds.
\eal
\ee
Therefore, as a bounded Leray-Hopf solution to \eqref{ueq} and $Z$ is smooth when $t<1$, $u$ is smooth when $t<1$ and bounded up to $t=1$. From \eqref{intequ}, $u$ is a bounded, classical solution to \eqref{ueq} in $\mathbb{R}^3 \times [0, 1)$. Recall that $Z$ becomes unbounded as $t \to 1^-$.
By \eqref{vzu}, $v=Z+u$ is the stated blow up solution to the forced Navier Stokes equation in $\mathbb{R}^3 \times [0, 1)$.
\qed
\medskip

Finally we remark that the constant $\delta_0$ can be estimated since it only depends on $T$,  the constants in the bound of the Stokes kernel and constants coming out of integrals of explicit functions.

\section*{Acknowledgments} We wish to thank Professors  Zijin Li,
  Xinghong Pan, Vladimir Sverak, Xin Yang and Na Zhao  for helpful discussions. The support of Simons
  Foundation grant 710364 is gratefully acknowledged.

\noindent {\it Statements and Declarations.} 1. There are no competing interests. 2. There is no data involved in this paper.

\bibliographystyle{plain}

%\bibliography{/Users/dykim/Dropbox/01Research/oralpaper}

%%%%%%%%%%%%%%%%%%%%%%%%%%%%%%%%%%%%
%\begin{comment}

\def\cprime{$'$}

\end{document}